\documentclass[12pt]{article}
\usepackage{amssymb}
\usepackage[leqno]{amsmath}
\usepackage{amsthm}
\usepackage{graphicx}
\usepackage{xr-hyper}
\usepackage[bookmarks,pdfnewwindow]{hyperref}

\hypersetup{
colorlinks=true,
linkcolor=black,
citecolor=black,
pdfauthor={Victor Ivrii},
pdftitle={Sharp Spectral Asymptotics for Magnetic Schr\"odinger Operator with Irregular Potential},
pdfsubject={Sharp Spectral Asymptotics},
pdfkeywords={Microlocal Analysis, Magnetic Schr\"odinger Operator},
baseurl={http://www.math.toronto.edu/ivrii/Research/preprints/},
}

\newcommand\W{{\rm W}}
\newcommand\MW{{\rm {MW}}}
\newcommand\Def{{\overset {\rm {def}}{\ =\ }}}
\newcommand\const{{\rm {const}}}

\newcommand\rank{\operatorname{rank}}

\newcommand\bR{{\mathbb R}}

\newcommand\bZ{{\mathbb Z}}
\newcommand\cA{{\mathcal A}}
\newcommand\cE{{\mathcal E}}
\newcommand\cC{{\mathcal C}}

\newcommand\cR{{\mathcal R}}
\newcommand\cZ{{\mathcal Z}}

\newtheorem{theorem}{Theorem}[section]

\theoremstyle{definition}

\newtheorem{remark}[theorem]{Remark}

\newcommand{\sect}[1]{\setcounter{equation}{0}\section{#1}}

\begin{document}

\title{%
\hfill{\normalsize\sl To Leonid Volevich in occasion of his 70th Birthday}\\
\bigskip
Sharp Spectral Asymptotics for Magnetic Schr\"odinger Operator with Irregular Potential.
}

\author{%
Victor Ivrii\footnote{Work was partially supported by NSERC grant OGP0138277.}
\footnote{Work was partially supported by Canada Council for the Arts via Killam Program.}
\date{August 2, 2004}}
\maketitle

{\abstract
In this paper I consider sharp spectral asymptotics for multidimensional magnetic Schr\"odinger operator with irregular coefficients with respect to two parameters -- semiclassical parameter $h$ and coupling parameter $\mu$. There are few principally different cases, depending on dimension, rank of magnetic intensity matrix, relation between $h$ and $\mu$ and some extra assumptions.
\endabstract}

\sect{Introduction}

\subsection{Preface}

In this paper I consider multidimensional Schr\"odinger operator
\begin{multline}
A=A_0+V(x),\qquad A_0=\sum_{j,k\le d}P_jg ^{jk}(x)P_k, \\
P_j=hD_j- \mu V_j(x),
\quad h\in (0,1],\ \mu \ge 1.
\label{0-1}
\end{multline}
It is characterized by {\sl magnetic field intensity matrix\/}
$(F_{jk})$ with $F_{jk}= (\partial _{x_j}V_k- \partial_{x_k}V_j)$, which is skew-symmetric $d\times d$-matrix, and $(F^j_m)=(g^{jk})(F_{km})$ which is unitarily equivalent to skew-symmetric matrix
$(g^{jk})^{\frac 1 2}(F_{jk})(g^{jk})^{\frac 1 2}$. Then all eigenvalues of
$(F^j_k)$ (with multiplicities) are $\pm i f_m$ ($f_m>0$, $m=1,\dots, r$) and $0$ of multiplicity $q=d-2r$ where $2r=\rank (F^j_k)$.

I formulate results only in the case of $g^{jk}=\delta_{jk}$ (thus covering the case $g^{jk}=\const$ as well) and $F_{jk}=\const$ (then one can select   linear vector potential  $(V_1(x),\dots,V_d(x))$). The results in the general case differ even in their statements. So we consider operator
\begin{equation}
\sum_{j,k}P_j^2+V(x), \qquad P_j=hD_j-\mu V_j(x)
\label{0-2}
\end{equation}
with linear functions $V_j(x)$ and constant magnetic intensity matrix $F_{jk}$.

The main goal of this paper is to present {\sl the local spectral asymptotics\/}, i.e. asymptotics  of 
\begin{equation}
\int e(x,x,\tau)\psi(x)\,dx
\label{0-3}
\end{equation}
as $h\to +0$, $\mu \to +\infty$ where $e(x,y,\tau)=e_{h,\mu}(x,y,\tau)$ is the Schwartz kernel of the spectral projector  of $A$ ($A$  is assumed to be a self-adjoint operator) and $\psi$  is a smooth cut-off function supported in the ball  $B(0,{\frac 1 2})$ (all the conjectures save self-adjointness are made for $B(0,1)$). Combined with  partition-rescaling arguments such asymptotics imply many asymptotics of eigenvalue counting function given by  formula  (\ref{0-3}) with $\psi=1$.

Under above assumptions in the appropriate coordinates after gauge transform \begin{equation}
V_j(x)=\left\{
\begin{aligned}
&-f_{j-r}x_{j-r}\qquad&&{\rm as\ }j=r+1,\dots, 2r\\
&0 &&{\rm otherwise}
\end{aligned}
\right.
\label{0-4}
\end{equation}
where $f_j>0$ and $\pm if_j$ are eigenvalues of $F^j_k$, $\rank F = 2r$.  Let $q=d-2r$.

In the special case $V=\const$ everything becomes explicit:  after $h$-Fourier transform with respect to $(x_{r+1},\dots,x_d)$ and 
change of variables $x_j\mapsto x_j+\mu^{-1}f_j^{-1}\xi_{2j+1}$ $(j=1,\dots,r)$  operator $A$ is transformed into
\begin{equation}
\sum_{1\le j \le r} \bigl(h^2D_j^2 + \mu^2 f^2_j x_j^2\bigr)+
\sum_{1\le k\le q}\xi_{2r+k}^2+V.
\label{0-5}
\end{equation}
Using Hermite function decomposition one can calculate easily $e(x,y,\tau)$ and prove that
\begin{multline}
e(x,x,\tau)= \cE_{d,r}^\MW (\tau)\Def \\
\omega_q(2\pi )^{-q} 
\sum _{\alpha \in {\bf Z}^{+r}}
\Bigl(\tau - \sum_j (2\alpha_j +1) f_j\mu h -V\Bigr)_+^{\frac q 2} 
\times 
\mu ^rh^{-d+r} f_1\cdots f_r\sqrt g
\label{0-6}
\end{multline}
where $\omega_q$ is a volume of the unit ball in $\bR^q$, $g=\det (g^{jk})^{-1}$.

In particular, for $q=0$ such operator has pure point spectrum of infinite multiplicity, consisting of {\sl Landau levels\/} 
$E_\alpha=\sum_j(2\alpha_j+1)f_j\mu h+V$ with $\alpha \in \bZ^{+\, r}$ while for $q\ge 1$ the spectrum is absolutely continuous and Landau levels are merely bottoms of its channels $[E_\alpha ,\infty)$.

Also, as $\mu h \gg 1$ one needs to include $-\sum_j f_j\mu h$ or even larger negative number in $V$ to avoid being in the classically forbidden zone.

Our goal is to study operator $A$ with variable potential $V(x)$; however the deep distinctions between 
cases $\rank F_{jk}=d$ and $\rank F_{jk}<d$ and between cases $\mu h\ll 1$, 
$\mu h\gg 1$ (and intermediate case $\mu h\sim 1$) are preserved as well as the importance of Landau levels. Also, $\cE^\MW$ gives a good approximation to $e(x,x,\tau)$ even if in some cases we need some corrections to derive sharper remainder estimate.

\subsection{Classical Dynamics}

Considering classical trajectories on some finite energy level one can prove easily that for $V=\const$ 

(i) As \ $\rank F=d=2$ \ particles move along circles of the radii 
$\asymp \mu^{-1}$.   As $\rank F=d=2r$, $r\ge 2$ trajectories are more complicated (generic trajectories are periodic if $f_1,\dots,f_r$ are commensurable or envelope tori otherwise) but they are confined to $C\mu^{-1}$-vicinities of their origins as well;

(ii) As $d=3$, $\rank F=2$ particles move along spirals. Similar description holds for $\rank F=2r < d$, $r\ge 1$: there is a cyclotronic movement as in (i) along $(x_1,\dots, x_{2r})$ and also a free movement along $x'''=(x_{2r+1},\dots,x_d)$.

So, we see a big difference between full- and not-full-rank cases and more subtle  difference between $d=2,3$ and $d\ge 4$.

As $V$ is variable 

(iii) As \ $\rank F=d=2$ \ §particles move approximately along circles of the radii 
$\asymp \mu^{-1}$ but the centers of these circles are drifting with the speed $\mu^{-1}|f_1|^{-1}\nabla V^\perp$ (thus the drift is orthogonal to the direction of the electric field). Similarly, as $\rank F=d=2r$ the trajectories described in (i) are drifting with the velocity $\mu^{-1}F^{-1}\nabla V$ where $F=(F_{jk})$ is the intensity matrix;

(iv) As  $d=3$, $\rank F=2$ there is a fast cyclotronic movement along $(x_1,x_2)$, temperate movement described by 1-dimensional Hamiltonian 
$\xi_3^2 + V (x',x_3)+ E_{\rm magn}$ along $x_3$ where $E_{\rm magn}$ is an energy of the fast movement (which is constant in our assumptions) and a slow drift along all variables. Similar description holds for $\rank F=2r<d$,
$r\ge 1$.

So, if $\rank F=d$ we can follow trajectories until time $T_1=\epsilon \mu$. Then  assuming that trajectories are non-periodic which happens if 
\begin{align}
&|V|\ge \epsilon_0,\label{0-7}\\
&|\nabla V|\ge \epsilon_0\label{0-8}
\end{align}
we hope to get remainder estimate which is $T_1^{-1}h=\mu^{-1}h$ times better than the principal part; since the latter is 
$\asymp h^{-d}\bigl(1+\mu^r h^r\bigr)$ 
we expect to derive asymptotics with $O\bigl(\mu^{-1}h^{1-d}(1+\mu^r h^r)\bigr)$ remainder estimate under assumptions (\ref{0-7})-(\ref{0-8}). Assuming sufficient smoothness (see below) it is a correct guess. 

However without condition (\ref{0-8}) trajectories could be periodic with period $T_0=\epsilon \mu^{-1}$ and we cannot expect remainder estimate  better than $O\bigl(\mu h^{1-d}\bigr)$ as $\mu h\le 1$.

On the other hand, if $\rank F<d$ we can follow trajectories until time $T_1=\epsilon$ and under assumptions (\ref{0-7})-(\ref{0-8}) these trajectories are non-periodic. So, we expect $O\bigl(h^{1-d}(1+\mu^r h^r)\bigr)$ remainder estimate then.

However, even without condition (\ref{0-8}) most  of the trajectories are non-periodic in a short run due to the free movement. Actually, in the classical settings periodic trajectories must have $\xi_{2r+1}=\dots =\xi_d=0$ and thus form a set of measure 0 but we are in semiclassics and talking about non-periodicity one should be able to observe the shift after ``quasi-period'' $\asymp \mu^{-1}$. Anyway, we guess that the lack of condition (\ref{0-8}) would not be that destructive (and even noticeable for not very large $\mu $). Both guesses are correct under some smoothness assumption.

As $\mu h\ll 1$ one can rescale to the {\sl standard case\/} (i.e. with $\mu=1$). Really, after scaling $\mu x\mapsto x$ and thus $hD\mapsto \mu h D$,
$\mu \mapsto 1$, $h\mapsto \mu h$ we get a standard case and scaling back we get the Weyl principal part 
\begin{equation}
\cE^\W(x,x,\tau) =\omega_d (\tau-V)_+^{-{\frac d 2}}h^{-d}\sqrt g
\label{0-9}
\end{equation}
and the remainder estimate $O(\mu h^{1-d})$. Much more useful however is a rescaling applied to the intermediate rather than the final results.

\subsection{Canonical form}

Starting again from appropriate coordinates and applying $h$-Fourier transform with respect to $x''=(x_{r+1},\dots,x_r)$ and change of variables 
$x_j\mapsto x_j+\mu^{-1}f_j^{-1}\xi_{2j+1}$ with respect to $x'=(x_1,\dots,x_r)$ we arrive instead of (\ref{0-5}) to 
\begin{multline}
\sum_{1\le j \le r} \bigl(h^2D_j^2 + \mu^2 f^2_j x_j^2\bigr)+
\sum_{1\le k\le q}(hD_{2r+k})^2+\\
V(x'+\mu^{-1}hF^{-1}D'', x''-\mu^{-1}hF^{-1}D',x''')
\label{0-10}
\end{multline}
where the third term is  $\mu^{-1}h$-pseudo-differential operator.

Note that on bounded energy levels operators $hD'$ and $\mu x'$ are now bounded. Then one can apply Taylor decomposition to the last operator in (\ref{0-10}) leading to an operator with the main part
\begin{equation}
\cA^0=\sum_{1\le j \le r} \bigl(h^2D_j^2 + \mu^2 f^2_j x_j^2\bigr)+
\sum_{1\le k\le q}(hD_{2r+k})^2+
V(\mu^{-1}hF^{-1}D'', x'',x'''),
\label{0-11}
\end{equation}
junior terms 
\begin{equation}
\cA'=\sum_{\alpha,\beta\in \bZ^{+\,r}, 1\le|\alpha|+|\beta|< l}
W_{\alpha\beta}(\mu^{-1}hF^{-1}D'', x'',x''') \times (\mu^{-1}h)^{|\beta|} x^{\prime\,\alpha}D^{\prime\,\beta}
\label{0-12}
\end{equation}
and the remainder estimate $O(\mu^{-l})$ (subject to the smoothness assumptions).

Applying Hermite decomposition with respect to $x'$ we see that as $q=0$ the main part  $\cA^0$ becomes a family of $\mu^{-1}h$-pseudododifferential operators with respect to $x''$
\begin{equation}
\cA_\alpha =V_\alpha=
\sum_{1\le k\le r} (2\alpha_j+1)f_j\mu h +
V(\mu^{-1}hF^{-1}D'', x'',x'''),
\label{0-13}
\end{equation}
while as $q\ge1$ it becomes a family of $q$-dimensional Schr\"odinger operators with respect to $x'''$\/ \ with potentials which are $\mu^{-1}h$-pseudododifferential operators:
\begin{align}
&\cA_\alpha =
\sum_{1\le k\le q}(hD_{2r+k})^2+V_\alpha, \label{0-14}\\
&V_\alpha=\sum_{1\le k\le rq} (2\alpha_j+1)f_j\mu h +
V(\mu^{-1}hF^{-1}D'', x'',x''').\notag
\end{align}

As $q=0$ the principal part of the spectral asymptotics for an individual operator is $\asymp \mu^r h^{-r}$ with the remainder estimate 
$O(\mu^{r-1} h^{1-r})$ under non-degeneracy condition
\begin{equation}
|V_\alpha|+|\nabla_{x'',\xi'''}V_\alpha|\ge \epsilon
\label{0-15}
\end{equation}
and $O(\mu^r h^{-r})$ otherwise (functions $V_\alpha $ could be very flat). Assume that this non-degeneracy condition (which is equivalent to (\ref{0-8}) as $\mu h\ll 1$ and replaces it otherwise)  holds.  Then there are 
$\asymp (\mu^{-r}h^{-r}+1)$ contributing operators in the family (the rest is excluded due to ellipticity arguments) thus leading to the principal part and the remainder estimate 
$\asymp h^{-d}$ and $O(\mu^{-1}h^{1-2r})$ respectively as $\mu h\le 1$ and 
$\asymp \mu^r h^{-r}$ and $O(\mu^{r-1}h^{-r+1})$ as  $\mu h\ge 1$.

For individual Schr\"odinger operator the principal part is 
$\asymp \mu^r h^{r-d}$ with the best possible remainder estimate 
$O(\mu^r h^{1+r-d})$ under condition (\ref{0-15}). However this condition is not crucial: this remainder estimate holds without such condition if either $q\ge 3$
or $q=2$, $V\in C^2$ and one can recover a weaker remainder estimate  otherwise.
Again there are $\asymp (\mu^{-r}h^{-r}+1)$ contributing operators in the family
thus leading to the principal part and the remainder estimate 
$\asymp h^{-d}$ and $O(h^{1-2r})$ respectively as $\mu h\le 1$ and 
$\asymp \mu^r h^{-r}$ and $O(\mu^rh^{-r+1})$ as  $\mu h\ge 1$ (in the best case).

Surely one can take care of junior terms (\ref{0-12}) and remove them if possible, thus reducing operator to Birkhoff normal form. It is easy in the smooth case as $d=2,3$; as $d=2$ the normal form is
\begin{equation}
\mu ^2 \sum_{m+j+k\ge 1} a_{mjk}(x_2, \mu ^{-1}hD_2) \cdot
 \bigl(x_1 ^2+\mu ^{-2} h ^2D_1 ^2\bigr)^m 
\mu ^{-2m-2j-k} h^k
\label{0-16}
\end{equation}
and for $d=3$ it is
\begin{multline}
\mu ^2 \sum_{m+p+j+k\ge 1} a_{mpjk}(x_2,x_3, \mu ^{-1}hD_2) \times\\
\times \bigl(x_1 ^2+\mu ^{-2} h ^2D_1 ^2\bigr)^m (hD_3)^{2p} 
\mu ^{-2m-2j-2p-k} h^k
\label{0-17}
\end{multline}
but as $r\ge 2$ the {\sl resonances} become one of the obstacles: namely we cannot remove terms 
\begin{align}
a_{\alpha,\beta}(x'',x''',\mu^{-1}hD'')\times &\prod_{1\le j\le r} (x_j+i\mu^{-1}hD_j)^{\alpha_j}(x_j-i\mu^{-1}hD_j)^{\beta_j}\times
(\xi''')^\gamma,\label{0-18}\\
\intertext{with}
&\sum_j(\alpha_j-\beta_j)f_j=0
\label{0-19}
\end{align}
($d=3$ is a  special case). So {\sl resonance\/} means that $\sum_j\gamma_jf_j=0$ with $\gamma\in \bZ^r$; $|\gamma|$ is {\sl an order of resonance}.

As $q=0$, the second order resonance terms are $(x_j+i\mu^{-1}\xi_j)(x_k-i\mu^{-1}\xi_k)$ with $f_j=f_k$, the third order resonance terms are $(x_j+i\mu^{-1}\xi_j)(x_k-i\mu^{-1}\xi_k)(x_m-i\mu^{-1}\xi_m)$
with $f_j=f_k+f_m$ and their conjugates, etc. As $q\ge1$ there are additional (but less malicious in the end of the day) resonance terms.

Another obstacle is the lack of a very large smoothness. I overcome both of these obstacles because most of the junior terms are non-essential small perturbations:
as $q\ge 2$ all junior terms are non-essential; as $q=1$ all terms but of the form 
$\mu^{-1}(x_j+i\mu^{-1}\xi_j)(x_k-i\mu^{-1}\xi_k)(x_m-i\mu^{-1}\xi_m)$ are are also non-essential (but such terms do not appear for special operator (\ref{0-2})), while for $q=0$ there are more essential terms.

Further, the canonical form reduction is possible for general operator (\ref{0-1}) as well but there will be essential cubic perturbations and third-order resonances become important as $q=1$.

\subsection{Main tools and classification} 

We will consider $u(x,y,t)$ the Schwartz kernel of operator $U(t)=e^{ih^{-1}t{\tilde A}}$ where ${\tilde A}$ is one of two framing approximations for $A$; in the case of irregular coefficients ${\tilde A}$ will include mollification as well. However
even in the case of smooth coefficients it is often convenient to take 
${\tilde A}$ different from $A$.

Microlocal analysis (propagation of singularities) of $u(x,y,t)$ in its original form or transformed (as $A$ is reduced to its canonical form) plays a crucial role.  The former approach is used in the case of {\sl weak\/}  magnetic field  and it is combined with  the standard theory rescaled.  In the latter approach is used in the case of {\sl intermediate\/} and stronger magnetic field and it is combined with the successive approximation method to construct solutions.

As $q=0$ and non-degeneracy condition holds we need to take mollification parameter $\varepsilon = C\mu h|\log h|$ to make a weak magnetic field approach working properly. On the other hand, if we reduce operator to its canonical form, we must take $\varepsilon = C(\mu^{-1}h|\log \mu|)^{\frac 1 2}$ as $d=2$ and $\varepsilon=C\max \bigl(\mu^{-1}, (\mu^{-1}h|\log \mu|)^{\frac 1 2}\bigr)$ as $d\ge 4$. Since we are interested in the smallest possible $\varepsilon$ we
apply weak magnetic field approach for $\mu \le C(h|\log h|)^{-{\frac 1 3}}$
as $d=2$ and for $\mu \le C(h|\log h|)^{-{\frac 1 2}}$ as $d\ge 2$. Otherwise
we refer to {\sl intermediate\/} magnetic field case as 
$\mu \le C(h|\log h|)^{-1}$, then {\sl strong\/} as $\mu \le \epsilon h^{-1}$,
{\sl superstrong\/} as $\mu \le Ch^{-1}$ and {\sl ultrastrong\/} as 
$\mu \ge Ch^{-1}$.

This classification holds as $q\ge 1$ but with a twist. If magnetic field is weak we study only ``original'' propagator $u(x,y,t)$ picking up 
$\varepsilon =C\rho^{-1} h|\log h|$
with $\rho =|\xi'''|$ as  $\rho \ge {\bar\rho}_1=C\max \bigl(\mu^{-1}, (\mu h|\log h|)^{\frac 1 2}\bigr)$ and $\varepsilon =C{\bar\rho}_1^{-1} h|\log h|$ otherwise. Inner zone $\{|\xi'''|\le {\bar\rho}_1\}$ is treated  by referring to the standard results rescaled.

If magnetic field is intermediate we apply a weak magnetic field approach in the outer zone $\{|\xi'''|\ge {\bar\rho}_1\}$. In the inner zone we apply reduction to the canonical form. Threshold between weak and intermediate magnetic field cases depends on $q$ and the presence of non-degeneracy condition. 

In the strong (and stronger) magnetic field cases we use  canonical form reduction everywhere and $\varepsilon = Ch|\log h|$.

Actually  the above choice of $\varepsilon$ is needed to provide a reduction; the further analysis as $q\ge 1$ requires to increase $\varepsilon$ and the actual choice of $\varepsilon$  varies (as $q=1,2$) depending on other assumptions and those indicated are kind of milestones rather than the actual values.

\sect{Weak Magnetic Field\/}

We assume that magnetic field is relatively weak and we prove our results by reducing to the standard case.

\subsection{Standard reference result rescaled} As I mentioned the intermediate result in the standard case will be valuable for us more than the final one; namely we have by rescaling 
\begin{equation}
|F_{t\to h^{-1}\tau}\bigl({\bar\chi}_T(t)\Gamma_x u \bigr)- 
\int T{\hat{\bar \chi}}\bigl((\tau -\lambda){\frac T h}\bigr)d_\lambda 
\cE^\MW (x,\lambda)|\le  Ch^{3-d}\quad\forall \tau:|\tau|\le\epsilon
\label{1-1}
\end{equation}
as $\mu \le \epsilon h^{-1}|\log h|^{-1}$, 
$Ch|\log h|\le T\le {\bar T}_0=\epsilon\mu^{-1}$;
here and below ${\bar\chi}$ is smooth and even satisfying admissibility condition of \cite{BrIvr} function, supported in $[-1,1]$ and equal 1 in $[-{\frac 1 2}, {\frac 1 2}]$,
${\hat{\bar \chi}}$ is its Fourier transform, $\Gamma_x v=v(x,x,t)$. 

We use also notations $\Gamma v =\int \Gamma_xv\, dx$ and 
$\phi_T(t)=\phi(t/T)$. Also we need to assume condition (\ref{0-7}). Estimate (\ref{1-1}) holds as coefficients of operator are either smooth or mollified  with $\varepsilon\ge Ch|\log h|$.

\subsection{Full-rank case}

We want to increase $T$ in (\ref{1-1}). To do this we analyze propagation of singularities. Let us introduce {\sl slow variables\/}
\begin{equation}
X_j=x_j -\sum \beta_{jk}p_k(x,\xi), \qquad p_k(x,\xi)=\mu^{-1}\xi_k-V_k 
\label{1-2}
\end{equation}
with $(\beta_{jk})$ inverse matrix to $(F_{jk})$.

Then the Poisson brackets satisfy
\begin{align}
&\{p_k,X_j\}=0\qquad &\forall j,k,\label{1-3}\\
&\{X_j,X_k\}=\mu^{-1}\beta_{jk},\  \{p_j,p_k\} = -\mu^{-1}hF_{jk}\qquad &\forall j,k.\label{1-4}
\end{align}

Then  in the classical dynamics for time $T$ the shift of 
$\phi (x,\xi)=\sum_j \ell_j X_j$ is equal to 
$\mu^{-1}T\bigl(\sum_{j,k}\beta_{jk}\ell_j\partial_kV+o(1)\bigr)$ as
$T\le {\bar T}_1=\epsilon \mu$  and under condition (\ref{0-8}) it will be of magnitude $\mu^{-1}T$ for an appropriate vector $\ell$, $|\ell|=1$. 

To make it observable from the point of view of microlocal analysis  one needs to satisfy {\sl the logarithmic uncertainty principle\/} (see \cite{BrIvr}) which in this case is
\begin{equation}
\mu^{-1}T \times \varepsilon \ge C\mu^{-1}h|\log h|
\label{1-5}
\end{equation}
because $X_1,\dots,X_d$ are linked by (\ref{1-5}) and semiclassical parameter is 
therefore $\mu^{-1}h$. We can avoid condition
\begin{equation}
\varepsilon \times \varepsilon \ge C\mu^{-1}h|\log h|
\label{1-6}
\end{equation}
because we do not use any reduction in this case; condition (\ref{1-6}) be very unpleasant for not very large $\mu$.

We want to satisfy (\ref{1-6}) with $T={\bar T}_0$ and therefore we pick up the
smallest $\varepsilon$ to satisfy this condition:
\begin{equation}
 \varepsilon = C\mu h|\log h|.
\label{1-7}
\end{equation}
Then we prove that
\begin{equation}
|F_{t\to h^{-1}\tau}\chi_T(t) \Gamma (u\psi) |
\le  Ch^s \qquad\forall \tau:|\tau|\le\epsilon\quad 
\forall T\in [{\bar T}_0,{\bar T}_1]
\label{1-8}
\end{equation}
where  $\psi=\psi(x)$ is a cut-off function described above, $\chi$ is smooth (and \cite{BrIvr} type) function supported in 
$[-1,-{\frac 1 2}]\cup[{\frac 1 2}, 1]$.

Then (\ref{1-1}) and (\ref{1-8}) imply that estimate
\begin{equation}
|F_{t\to h^{-1}\tau}{\bar\chi}_T(t) \Gamma(u\psi)  - \int\int
T{\hat{\bar \chi}}\bigl((\tau -\lambda){\frac T h}\bigr)d_\lambda \cE^\MW (x,\lambda)\,\psi(x)dx|\le  Ch^{3-d}
\label{1-9}
\end{equation}
holds with any $T\in [{\bar T}_0,{\bar T}_1]$ and 
$\forall \tau:|\tau|\le\epsilon$.

We must here and below integrate with respect to $\psi(x)\,dx$ because $X_j$ contain $\xi_k$ and we need a real trace $\Gamma$ and not just a restriction to the diagonal $\Gamma_x$.

Then Tauberian arguments imply spectral asymptotics with the remainder estimate 
$O({\bar T}_1^{-1}h^{1-d})=O(\mu^{-1}h^{1-d})$ while an approximation error is 
$O(\varepsilon^l |\log \varepsilon|^{-\sigma}h^{-d})$ provided
\begin{equation}
V\in C^{l,\sigma}\qquad (l,\sigma)\succeq (1,1)
\label{1-10}
\end{equation}
where $C^{l,\sigma}$ is a class of functions with $k$-th derivatives  continuous with continuity modulus $t^{l-k}|\log t|^{-\sigma}$ functions, 
$k=\lfloor l\rfloor$ unless $l\in \bZ$, $\sigma \le 0$ when $k=l-1$.

Thus we arrive to our first statement:

\begin{theorem}
\label{thm-1-1}
Let $d=2r$, $V\in C^{l,\sigma}$ with $(l,\sigma)\succeq (1,2)$ and conditions $(\ref{0-7}),(\ref{0-8})$  hold. Then for
\begin{equation}
\mu \le h^{\delta -1}
\label{1-11}
\end{equation}
with an arbitrarily small exponent $\delta >0$ the following estimate holds
\begin{multline}
\cR\Def |\int \bigl(e(x,x,0) - \cE^\MW (x,0)\bigr)\psi(x)\,dx |\le\\
C\mu^{-1}h^{1-d}+ C(\mu h|\log h|)^l |\log h|^{-\sigma} h^{-d}.
\label{1-12}
\end{multline}

In particular, as   $(l,\sigma)=(3,1)$  sharp remainder estimate 
$\cR\le C \mu^{-1}h^{1-d}$ holds for 
$\mu \le h^{-{\frac 1 2}}|\log h|^{-{\frac 1 2}}$
and as   $(l,\sigma)=(2,1)$  this sharp remainder estimate holds for 
$\mu \le h^{-{\frac 1 3}}|\log h|^{-{\frac 1 3}}$.
\end{theorem}

\subsection{Non-full-rank case}

As $q\ge 1$ we can use a free movement along $x'''$ to extend (\ref{1-9}) from $T={\bar T}_0$ to larger $T$. Let us introduce $\rho$-partition with respect to $\xi'''$ with  $\rho={\frac 1 2}|\xi'''|$ and let us consider a partition element with $\rho\ge C\mu^{-1}$. Then for time $T\le {\bar T}_1=\epsilon \rho$ the shift with respect to $x'''$ will be of magnitude $\epsilon \rho$ and in order to be observable it must satisfy logarithmic uncertainty principle
\begin{equation}
\rho \times \rho T\ge Ch|\log h|\iff T\ge C\rho^{-2}h|\log h|.
\label{1-13}
\end{equation}
In order to be able to plug $T={\bar T}_0=\epsilon \mu^{-1}$ we must take
\begin{equation}
\rho\ge {\bar\rho}_1=C\max\bigl(\mu^{-1},(\mu h|\log h|)^{\frac 1 2}\bigr).
\label{1-14}
\end{equation}
Then both (\ref{1-8}) and (\ref{1-9}) hold for $u \psi$ replaced by $Q(hD_x)(u\psi)$ with $h$-pseudodifferential oprator $Q$ with symbol supported in the indicated partition element. 

Then due to Tauberian arguments the contribution of this element to the remainder estimate is 
$Ch^{1-d}{\bar T}_1^{-1}\rho^q$ and therefore the total contribution of the 
{\sl outer\/} zone $\{|\xi'''|\ge {\bar\rho}_1\}$
to the remainder estimate does not exceed 
$ Ch^{1-d}\int {\bar T}_1^{-1}(\rho)\rho^{q-1}\,d\rho$ which is $O(h^{1-d})$ as $q\ge 2$ and $O(h^{1-d}\log \mu)$ as $q=1$. To get rid off the logarithmic factor we  prove that for $(l,\sigma)\succeq (1,2)$ one can take 
${\bar T}_1=\epsilon \rho|\log \rho|^2$ in the propagation.

Also in this zone $\varepsilon$ must satisfy logarithmic uncertainty principle 
$\varepsilon \times\rho \ge Ch|\log h|$ and taking 
$\varepsilon = C\rho^{-1}h|\log h|$ there we get mollification error
$Ch^{-d}\int \varepsilon ^l|\log \varepsilon|^{-\sigma}\rho^{q-1}\,d\rho$
which is $O(h^{1-d})$ as either $q\ge 2$, $(l,\sigma)\succeq (1,1)$ or 
$q=1, (l,\sigma)\succeq (1,2)$. Note that now mollification parameter $\varepsilon$ depends on $(x,\xi)$.

On the other hand, contribution of the zone $\{|\xi'''|\le {\bar\rho}_1\}$ to the remainder estimate is $C{\bar T}_0^{-1}{\bar\rho}_1^qh^{1-d}= 
O\bigl(h^{1-d}+(\mu h )^{{\frac q 2}+1}|\log h|^{\frac q 2}h^{-d}\bigr)$.

Finally, in this zone we pick $\varepsilon = C{\bar\rho}_1^{-1}h|\log h|$ and its contribution  to the mollification error will be less than what we already got. 

So we arrive to estimate
\begin{equation}
\cR \le Ch^{1-d}+C(\mu h )^{{\frac q 2}+1}|\log h|^{\frac q 2}h^{-d}
\label{1-15}
\end{equation}
which actually could be improved to

\begin{theorem}
\label{thm-1-2}
 Let condition $(\ref{0-7})$ be fulfilled.  Let   either $q= 1$, $V\in C^{1,2}$ or $q\ge 2$, $V\in C^{1,1}$. Then 
\begin{equation}
\cR\le 
Ch^{1-d}+ C(\mu h)^{{\frac q 2}+1}  h^{-d}.
\label{1-16}
\end{equation}

In particular, for  $\mu \le \mu^*_{1(q)} =h^{-{\frac q {q+2}}}$  sharp remainder estimate $\cR\le C h^{1-d}$ holds. As $q\ge 2$ this sharp remainder estimate holds for $\mu \le h^{-{\frac 1 2}}$.
\end{theorem}

Note, that at this stage extra smoothness is not very useful. Furthermore we can assume in what follows that $\mu \ge \mu^*_{1(q)}$ and thus 
${\bar\rho}_1=C(\mu h|\log h|)^{\frac 1 2}$.
 
To improve estimates (\ref{1-15}), (\ref{1-16}) one needs to use better arguments in  the  inner zone $\cZ=\{|\xi'''|\le {\bar\rho}_1\}$. Under non-degeneracy condition (\ref{0-8}) we can apply the same arguments as in the proof of theorem \ref{thm-1-1} as long as $\varepsilon = C\mu h|\log h|$ in $\cZ$. Then contribution of this zone to the remainder estimate becomes $O(h^{1-d})$ while its contribution to an approximation  error becomes 
$O({\bar\rho}_1^q \varepsilon ^l|\log \varepsilon|^{-\sigma}h^{-d})$ and we arrive to

\begin{theorem}
\label{thm-1-3}
 Let  condition $(\ref{0-7}), (\ref{0-8})$ be fulfilled.
 Let    $q\ge 1$, $V\in C^{l,\sigma}$ and either $q=1$, 
 $(l,\sigma)\succeq (1,2)$ or $q\ge 2$,  $(l,\sigma)\succeq (1,1)$ . Then  
\begin{equation}
\cR\le Ch^{1-d}+ C(\mu h)^{{\frac q 2}+l} |\log h|^{l-\sigma}  h^{-d}.
\label{1-17}
\end{equation}

In particular, sharp remainder estimate $\cR\le C h^{1-d}$ holds for  $\mu \le h^{-{\frac 1 2}}|\log h|^{-{\frac 1 2}}$
as  $q=1$, $(l,\sigma)=({\frac 3 2},{\frac 1 2})$.
 \end{theorem}

\sect{Intermediate and strong magnetic field}

Now magnetic field is stronger than before but still either below 
$\epsilon h^{-1}|\log h|^{-1}$ or between this value and $\epsilon h^{-1}$. There is certain difference in the analysis of these two cases and for general operator (\ref{0-1}) some statements would slightly differ as well.

\subsection{Full-rank case}
In this case we pick 
\begin{equation}
\varepsilon = \left\{\begin{aligned}
&C(\mu^{-1}h|\log h|)^{\frac 1 2}\qquad &r=1,\\
&C\max \bigl(\mu^{-1},(\mu^{-1}h|\log h|)^{\frac 1 2}\bigr)
&r\ge2
\end{aligned}\right.
\label{2-1}
\end{equation}
and reduce operator to a canonical form in the smooth case or ``a poor man's canonical form'' otherwise. Cases $r=1$ and $r\ge 2$ differ because of the reduction: we need to solve equation
\begin{equation}
\{a_0,S\} = V-W,\qquad a_0=\sum_{1\le j\le r}f_j(x_j^2+\xi_j^2)
\label{2-2}
\end{equation}
where $W$ consists of unremovable terms; as $r=1$ this equation  is solved by integration along circles leaving 
$W=W\bigl(x_2,\xi_2,(x_1^2+\xi_1^2)^{\frac 1 2}\bigr)$ while for $r\ge 2$ it is solved by Taylor decomposition and we must assume that 
$\varepsilon \ge C\mu^{-1}$. On the other hand, we always need to assume  (\ref{1-6}) now because we consider $x'',\xi''$ as dual variable and we need to consider $\mu^{-1}h$-pseudodifferential operators.

In the best possible case we would get something 
similar the family of separate scalar $\mu^{-1}h$-pseudodifferential operators but the same results hold in the general case as well:

\begin{theorem}
\label{thm-2-1}
Let $d=2r$, $V\in C^{l,\sigma}$ with $(l,\sigma)\succeq (1,2)$ and conditions $(\ref{0-7}),(\ref{0-8})$  hold. Then 

(i) For $\mu \le  h^{-1}|\log h|^{-1}$ estimate
\begin{equation}
\cR\le
C\mu^{-1}h^{1-d}+ C\varepsilon ^l |\log h|^{-\sigma} h^{-d}\qquad \forall \tau:|\tau|\le \epsilon
\label{2-3}
\end{equation}
holds with $\varepsilon=\mu^{-1}$; 

(ii) For 
$h^{-1}|\log h|^{-1}\le \mu \le \epsilon h^{-1}$ estimate $(\ref{2-3})$ holds with $\varepsilon=C(\mu^{-1}h|\log h|)^{\frac 1 2}$.

In particular, as   $(l,\sigma)=(3,1)$  sharp remainder estimate estimate 
$\cR\le C \mu^{-1}h^{1-d}$ holds for 
$h^{-{\frac 1 2}}|\log h|^{-{\frac 1 2}}\le \mu \le \epsilon h^{-1}$;

(iii) As $r=1$,  estimate $(\ref{2-3})$ holds with 
$\varepsilon=(\mu^{-1}h|\log h|)^{\frac 1 2}$ but for $\cR_I$ instead of 
$\cR$; here and below $\cR_I$ is defined by formula $(\ref{1-12})$ but with $\cE^\MW$ replaced
by $\cE^\MW_I$ which is defined by $(\ref{0-6})$ with $V(x)$ replaced by $W(x)$
where $W(x)$ is an average of $V(y)$ along circle 
$\cC_x=\{y:|x-y|=(\mu f)^{-1}(\tau -V(x))^{\frac 1 2}\}$; $f$ is a scalar intensity of magnetic field.

In particular, as   $(l,\sigma)=(2,1)$  sharp remainder estimate estimate 
$\cR_I\le C \mu^{-1}h^{1-d}$ holds for 
$\mu \ge h^{-{\frac 1 3}}|\log h|^{-{\frac 1 3}}$.
\end{theorem}

This statement together with theorem \ref{thm-1-1} cover case
$\mu \le \epsilon h^{-1}$ completely.

\subsection{Non-full rank case. I} 
The same classification and definition of $\varepsilon$ persist as $q\ge 1$; however we reduce operator to a canonical form only in inner zone 
$\cZ=\{|\xi'''|\le {\bar\rho}_1\}$ (in the intermediate magnetic field case); in 
the outer zone we apply the weak magnetic field approach and pick up
$\varepsilon = C|\xi'''|^{-1}h|\log h|$; furthermore after reduction is done,
$\varepsilon$ is redefined (increased) in the inner zone as well. 

There are few different statements to prove; the first one is a generic one:

\begin{theorem}
\label{thm-2-2}
Let $q\ge 1$, $V\in C^{l,\sigma}$. Then 

(i) As $q\ge 3$, $(l,\sigma)=(1,1)$ sharp remainder estimate $\cR\le Ch^{1-d}$ holds for $\mu \le Ch^{-1}$;

(ii) As $q=2$, $(l,\sigma)=(1,1)$ remainder estimate 
\begin{equation}
\cR\le Ch^{1-d}+ C\mu h^{{\frac 5 3}-d}
\label{2-4}
\end{equation} 
holds for $\mu \le Ch^{-1}$;

(iii) As $q=1$, $(l,\sigma)\succeq (1,2)$ remainder estimate
\begin{equation}
\cR\le  C\mu h^{{\frac 4 3}-d}+
C\mu ^{1-{\frac l 2}}h^{1-d}|\log h|^{-{\frac \sigma 2}}
\label{2-5}
\end{equation} 
holds for $h^{-{\frac 1 2}}|\log h|^{-{\frac 1 2}}\le \mu \le Ch^{-1}$;

(iv) As $d=3$, $(l,\sigma)=(1,2)$ remainder estimate
\begin{equation}
\cR_I\le  C\mu h^{{\frac 4 3}-d}
\label{2-6}
\end{equation} 
holds for $h^{-{\frac 1 3}}\le \mu \le Ch^{-1}$; here $\cR_I$ is again defined
by $(\ref{1-12})$ with $\cE^\MW$ replaced by $\cE^\MW_I$ which is defined by $(\ref{0-6})$ with $V(x)$ replaced by $W(x)$
where $W(x)$ is an average of $V(y)$ along circle 
$\cC_x=\{y:x_3=y_3, |x-y|=(\mu f)^{-1}(\tau -V(x))^{\frac 1 2}\}$ as magnetic field is directed along $x_3$.
\end{theorem}

In what follows we need to treat only cases $q=1,2$. As non-degeneracy condition is fulfilled we get

\begin{theorem}
\label{thm-2-3}
Let $q=1,2$, $V\in C^{l,\sigma}$ and conditions $(\ref{0-7}),(\ref{0-8})$  hold. Then  

(i) As $q=2$, $(l,\sigma)=(1,1)$ sharp remainder estimate $\cR\le Ch^{1-d}$ holds for $\mu \le Ch^{-1}$;

(ii) As $q=1$, $(l,\sigma)\succeq (1,2)$ remainder estimate
\begin{equation}
\cR\le  C h^{1-d}+
C\mu ^{{\frac 1 2}-l}h^{{\frac 1 2}-d}|\log h|^{{\frac 1 2}- \sigma }
\label{2-7}
\end{equation} 
holds for $h^{-{\frac 1 2}}|\log h|^{-{\frac 1 2}}\le \mu \le Ch^{-1}$;

(iii) As $d=3$, $(l,\sigma)=(1,2)$ sharp remainder estimate
$\cR_I\le Ch^{1-d}$ holds for $\mu \le Ch^{-1}$.
\end{theorem}

So far all the the results of the article could be generalized to a general operator (\ref{0-1}) with modification of non-degeneracy condition (\ref{0-8})
and, as $q=0,1$, $r\ge 2$ with the special attention to the third order resonances (because in the general case they could lead to non-removable $O(\mu^{-1})$ terms in the canonical form).

\subsection{Non-full-rank case. II}

In this section we exploit more specific properties of operator (\ref{0-2}), namely that $f_j$ have constant multiplicities.

\begin{theorem}
\label{thm-2-4}
Let $q=1,2$, $V\in C^{l,\sigma}$ and condition $(\ref{0-7})$  hold. Then  

(i) As $q=2$, $(l,\sigma)\succeq (1,1)$  remainder estimate 
\begin{equation}
\cR\le Ch^{1-d}+ 
C\mu h^{{\frac {2l}{l+2}}+1-d} |\log h|^{-{\frac {2\sigma}{l+2}}}
\label{2-8}
\end{equation} 
holds for $\mu \le Ch^{-1}$; in particular sharp remainder
estimate $\cR\le Ch^{1-d}$ holds as $(l,\sigma)=(2,0)$.

(ii) As $q=1$, $(1,2)\preceq (l,\sigma)\preceq (2,0)$ remainder estimate
\begin{equation}
\cR\le Ch^{1-d}+ 
C\mu h^{{\frac {l}{l+2}}+1-d} |\log h|^{-{\frac {\sigma}{l+2}}}+
C\mu ^{1-{\frac l 2}}h^{1-d}|\log h|^{-{\frac \sigma 2}}
\label{2-9}
\end{equation}
holds for $h^{-{\frac 1 2}}|\log h|^{-{\frac 1 2}}\le \mu \le Ch^{-1}$;

(iii) As $d=3$, $(1,2)\preceq (l,\sigma)\preceq (2,0)$ remainder estimate
\begin{equation}
\cR_I\le Ch^{1-d}+ 
C\mu h^{{\frac {l}{l+2}}+1-d} |\log h|^{-{\frac {\sigma}{l+2}}}
\label{2-10}
\end{equation}
holds for $h^{-{\frac 1 3}}|\log h|^{-{\frac 1 3}}\le \mu \le Ch^{-1}$.
\end{theorem}

Further, as $r\ge 2$ Diophantine properties of $(f_1,\dots,f_r)$ can play role.
Assume that 
\begin{equation}
{\bf n}( \hbar,\tau)\Def\#\bigl\{\alpha \in \bZ^{+\, r}, 
\sum_j (2\alpha _j+1)f_j +V\hbar < \tau\bigr\}
\label{2-11}
\end{equation}
satisfies estimate
\begin{multline}
|{\bf n}( \hbar,\tau)-{\bf n}( \hbar,\tau')|\le 
C\hbar^{-r}\bigl(|\tau -\tau'|+\nu (\hbar)\bigr) \label{2-12}\\
\forall \hbar\in(0,1]\;\forall \tau,\tau':|\tau|\le C, |\tau'|\le C
\end{multline}
with $\nu(\hbar)=o(\hbar)$ (it holds with $\nu (\hbar)=\hbar $ for sure). Two following theorems improve theorems \ref{thm-2-2}, \ref{thm-2-4} respectively:

\begin{theorem}
\label{thm-2-5}
Let $r\ge 2$, $q=1,2$, $V\in C^{l,\sigma}$ and conditions $(\ref{0-7}),(\ref{2-11})$  hold. 
Then

(i) As $q=2$, $(l,\sigma)\succeq (1,1)$  remainder estimate 
\begin{equation}
\cR\le Ch^{1-d}+ C\nu (\mu h) h^{{\frac 2 3}-d}
\label{2-13}
\end{equation} 
holds for $\mu \le Ch^{-1}$; 

(ii) As $q=1$, $(l,\sigma)\succeq (1,2)$  remainder estimate 
\begin{align}
\cR\le &Ch^{1-d}+ 
C\mu ^{{\frac 1 2}-l}h^{{\frac 1 2}-d}|\log h|^{{\frac 1 2}- \sigma }+ \label{2-14}\\
&C\nu (\mu h)  \Bigl( h^{{\frac 1 3}-d}+
\mu ^{-{\frac l 2}}h^{-d}|\log h|^{-{\frac \sigma 2}}\Bigr)
\notag
\end{align} 
holds for $\mu \le Ch^{-1}$. 
\end{theorem}

\begin{theorem}
\label{thm-2-6}
Let $r\ge 2$, $q=1,2$, $V\in C^{l,\sigma}$ and conditions $(\ref{0-7}),(\ref{2-11})$  hold. 
Then  

(i) As $q=2$, $(l,\sigma)\succeq (1,1)$  remainder estimate 
\begin{equation}
\cR\le Ch^{1-d}+ 
C\nu (\mu h)  h^{{\frac {2l}{l+2}}-d} |\log h|^{-{\frac {2\sigma}{l+2}}}
\label{2-15}
\end{equation} 
holds for $\mu \le Ch^{-1}$; 

(ii) As $q=1$, $(l,\sigma)\succeq (1,2)$  remainder estimate 
\begin{align}
\cR\le &Ch^{1-d}+  
C\mu ^{{\frac 1 2}-l}h^{{\frac 1 2}-d}|\log h|^{{\frac 1 2}- \sigma }+
\label{2-16}\\
&C\nu (\mu h)  \Bigl(h^{{\frac {l}{l+2}}-d} |\log h|^{-{\frac {\sigma}{l+2}}}+
\mu ^{-{\frac l 2}}h^{-d}|\log h|^{-{\frac \sigma 2}}\Bigr)
\notag
\end{align} 
holds for $\mu \le Ch^{-1}$. 
\end{theorem}

Note that in estimates (\ref{2-13})--(\ref{2-16}) right-hand expressions are sums of the right-hand expressions under non-degeneracy condition and of the right-hand expressions without microhyperbolicity conditions, but latter are multiplied by $\nu (\mu h) (\mu h)^{-1}$.

\sect{Superstrong and ultrastrong magnetic field}

In this case $\mu \ge \epsilon h^{-1}$ and the distance Landau levels increases; in the case of the ultrastrong magnetic field only one (may be multiple) level should be considered.

\subsection{Superstrong magnetic field}
In this case $h^{-1}\epsilon \le \mu  \le Ch^{-1}$ and the magnitude of the principal part of the asymptotics is still the same ($h^{-d}$) as well as the remainder estimates in theorems \ref{thm-2-1}--\ref{2-4}. The only difference is that non-degeneracy condition (\ref{0-8}) is replaced by 
\begin{equation}
|\tau -V-\sum_j (2\alpha_j+1)f_j\mu h |+|\nabla V|\ge \epsilon_0\qquad \forall \alpha \in \bZ^{+\,r}.
\label{3-1}
\end{equation}

\begin{theorem}
\label{thm-3-1}

(i) Statements of theorems \ref{thm-2-2}, \ref{thm-2-4} remain true for
$h^{-1}\epsilon \le \mu  \le Ch^{-1}$;

(ii) Statements of theorems \ref{thm-2-1}, \ref{thm-2-3} remain true for
$h^{-1}\epsilon \le \mu  \le Ch^{-1}$ with condition $(\ref{0-8})$ replaced by
$(\ref{3-1})$ with $\tau=0$ and condition $(\ref{0-7})$ skipped;

(iii) As $q=0$ under condition 
\begin{equation}
|\tau -V-\sum_j (2\alpha_j+1)f_j\mu h  |\ge  \epsilon_0 \qquad \forall \alpha \in \bZ^{+\,r}
\label{3-2}
\end{equation}
with $\tau=0$ estimate $\cR\le Ch^s$ holds with arbitrarily large $s$ (spectral gaps).
\end{theorem}

\begin{remark}
\label{rem-3-2}
As $1\le \mu \le Ch^{-1}$ under condition 
\begin{equation}
\tau -V-\sum_j  f_j \mu h\le  -\epsilon_0 
\label{3-3}
\end{equation}
estimate $\cR \le Ch^s$ holds with $\cE^\MW (x,\tau)=0$ ($\tau =0$) and with arbitrarily large  $s$.
\end{remark}

\subsection{Ultrastrong magnetic field}

In this case $\mu\ge Ch^{-1}$ and in order not to be below the bottom of the spectrum one should modify condition $V\in C^{l,\sigma}$. Assume instead that
\begin{align}
&V=-\sum_j (2{\bar\alpha}_j+1)f_j \mu h + W,\quad &{\bar\alpha}\in \bZ^{+\,r},\ W\in C^{l,\sigma}  \qquad (q=0);
\label{3-4}\\
&V=-\sum_j f_j \mu h + W,\quad &W\in C^{l,\sigma}\qquad (q\ge 1).
\label{3-5}
\end{align}

\begin{theorem}
\label{thm-3-3} Let $q=0$ and conditions $(\ref{3-4})$, $(\ref{3-1})$ be fulfilled with $\tau=0$. Then

(i) Estimate 
\begin{equation}
\cR\le C\mu^{r-1}h^{-r+1};
\label{3-6}
\end{equation}
holds for $\mu \ge Ch^{-1}$;

(ii) Furthermore, under condition $(\ref{3-2})$ estimate 
$\cR\le C\mu^{-s}$ holds with arbitrarily large $s$.
\end{theorem}

\begin{theorem}
\label{thm-3-4} Let $q\ge 1$ and condition $(\ref{3-5})$ be fulfilled. Then

(i) Estimate 
\begin{equation}
\cR\le C\mu^r h^{-r+1}
\label{3-7}
\end{equation}
holds for $\mu \ge Ch^{-1}$, $q\ge 3$;

(ii) Estimate 
\begin{equation}
\cR\le C\mu^r h^{-r+1}\bigl(1 + h^{-1+{\frac {lq} {l+2}}}|\log h|^{-{\frac {\sigma q}{l+2}}}\bigr)
\label{3-8}
\end{equation}
holds for $\mu \ge Ch^{-1}$, $q=1,2$; in particular, as $q=2$, $(l,\sigma)=(2,0)$ sharp remainder estimate $(\ref{3-7})$ holds.

(iii) Under condition $(\ref{3-1})$  sharp remainder 
estimate $(\ref{3-7})$ holds for $\mu \ge Ch^{-1}$, $q=1,2$.

(iv) Under condition $(\ref{3-3})$ estimate $\cR\le C\mu^r h^s$ holds with $\cE^\MW (x,\tau)=0$ ($\tau =0$) and with arbitrarily large $s$.

\end{theorem}

\sect{Remarks. Generalizations}

One can generalize the results stated above.

\begin{remark}
\label{rem-4-1}
One can get rid of condition (\ref{0-7}) by method of rescaling; then all the results remain the same.
\end{remark}

\begin{remark}
\label{rem-4-2}
Instead of operator (\ref{0-2}) one can consider operator (\ref{0-1}) in its full generality:

(i) as $q\ge 3$ no modification in conditions is needed.

(ii)  As $q=0,1,2$ non-degeneracy conditions (\ref{0-8}) and (\ref{3-1}) should be modified; as $r=1$ one needs to replace $|\nabla V|$ by 
$|\nabla (V/f)|$ with $f=f_1$ ($\tau=0$); as $r\ge 2$ the modification is more profound because we are essentially in the matrix situation. For example  if $f_1,\dots,f_d$
have constant multiplicities our condition looks like
\begin{equation}
|  V + \sum_j f_j\tau_j |+ 
|\nabla \bigl(V^{-1}( \sum_j f_j\tau_j)\bigr)|\ge \epsilon_0\qquad \forall \tau_1\ge0,\dots,\tau_r\ge 0
\label{4-1}
\end{equation}
as ($\tau=0$) and we need to assume that (\ref{0-7}) holds.

(iii) Smoothness conditions to $g^{jk}$, $F_{jk}$ should be at least $(l,\sigma)$ and also at least $(2,1)$ (required for reduction arguments) but for they could be even stronger to get  a proper mollification error.

(iv) Constant multiplicity of $f_j$ which was taken for granted is no more guaranteed and the results of subsection 2.3 require it.

(v) Further, condition (\ref{2-12}) should be fulfilled for 
${\bf n}(x,\hbar)$ integrated with respect to $x$ and it is fulfilled automatically 
with $\nu(\hbar)=\hbar^\kappa$ provided $\{\nabla (f_2/f_1),\dots, \nabla (f_r/f_1)\}$ has at least rank $\kappa -1$.

(vi) As $q=1$ and non-degeneracy condition is not fulfilled, unremovable 
$O(\mu^{-1}))$ terms in the canonical form could lead to a some correction term
in order to save the estimate.
These terms can appear due to variable $g^{jk}$, $F_{jk}$ even if $f_j=\const$ and they are unremovable due to the third-order resonances.

(vii) In the case of $q=0$ and the ultrastrong magnetic field one should require that for each $\alpha\in \bZ^{+\,r}$ such that
$|\sum_j ({\bar\alpha}_j-\alpha_j)f_j| \le \epsilon$  \;
for each $j$ either
${\bar\alpha}_j=\alpha_j$ or $f_j$ has a constant multiplicity.

\end{remark}

\begin{remark}
\label{rem-4-3}
(i) The results of this article are proven in three papers \cite{IRO3}, \cite{IRO4}, \cite{IRO5} where the first one is dealing with 2,3-dimensional cases and the second and the third are dealing with the higher dimensions.  One can access  
them from

\centerline{\tt http://www.math.toronto.edu:/ivrii/Research/preprints.html}

\noindent
as well as relevant talks to be viewed on a computer screen rather than printed.

(ii) I my forthcoming papers I am planning to get rid off assumption 
``$\rank (F_{jk})=\const$''; surely results of the case $q=0$ will be no more valid.
\end{remark}

\bibliographystyle{amsalpha}

\providecommand{\bysame}{\leavevmode\hbox to3em{\hrulefill}\thinspace}

\vglue .06truein


\vglue .06truein

\begin{tabular}{rrl}
&{\hskip 220 pt} &Department of Mathematics,\cr
&&University of Toronto,\cr
&&100  St.George Str.,\cr
&&Toronto, Ontario M5S 3G3\cr
&&Canada\cr
&&ivrii@math.toronto.edu\cr
&&Fax: (416)978-4107\cr
\end{tabular}

\end{document}